\documentclass[12pt]{amsart}
\usepackage[dvips]{graphicx}
\usepackage{epsfig}
\usepackage{amssymb}
\usepackage[usenames, dvipsnames]{color}
\usepackage{hyperref}
\usepackage{verbatim}
\usepackage[normalem]{ulem}

\theoremstyle{plain}
\newtheorem{theorem}{Theorem}[section]

\theoremstyle{definition}

\theoremstyle{remark}

\newcommand{\bR}{{\mathbb R}}

\def\th{\theta}

\def\p{\partial}

\def\na{\nabla}
\def\al{\alpha}

\def\O{\Omega}

\def\be{\begin{equation}}
\def\ee{\end{equation}}
\def\bes{\begin{equation*}}
\def\ees{\end{equation*}}
\def\bali{\begin{aligned}}
\def\eali{\end{aligned}}
\def\al{\begin{aligned}}
\def\eal{\end{aligned}}

\def\lab{\label}

\def\2O{\underline{\O}}

\numberwithin{equation}{section}

\makeatletter
\def\dashint{\operatorname%
{\,\,\text{\bf--}\kern-.98em\DOTSI\intop\ilimits@\!\!}}
\makeatother

\begin{document}


\title[blow up solution]{A blow up solution of the Navier-Stokes equations with a super critical forcing term}

\thanks{}

\author[Q. S. Zhang]{Qi S. Zhang}

\address[]{Department of mathematics, University of California, Riverside, CA 92521, USA}

\email{qizhang@math.ucr.edu}



\begin{abstract}
 A forced solution $v$ of the axially symmetric Navier-Stokes equation in a finite cylinder $D$ with suitable boundary condition is constructed. The forcing term, whose order of scaling is slightly worse than the critical order $-2$,  is in the mildly super critical space $L^q_t L^1_x$ for all $q>1$. The velocity, which is smooth until its final blow up moment, is in the energy space through out.

\end{abstract}
\maketitle


\section{Statement of result and Proof}

The  Navier-Stokes equations (NS) describing the motion of viscous incompressible fluids in a domain $D \subset \bR^3$ can be written as
\be
\lab{nse}
\mu \Delta v -  v \nabla v - \nabla P -\partial_t v =f, \quad div \, v=0, \quad \text{in} \quad D \times (0, \infty),  \quad v(\cdot, 0)=v_0(\cdot)
\ee The unknowns are the velocity $v$ and the pressure $P$.
 $\mu>0$ is the viscosity constant, which will be taken as $1$; $f$ is a given forcing term. In order to solve the equation,  suitable boundary conditions should be given if $\partial D \neq \emptyset$.

 Thanks to Leray \cite{Le2}, if $D=\bR^3$, $v_0 \in L^2(\bR^3)$, $f=0$,  the Cauchy problem has a solution in the energy space $(v, \nabla v) \in   (L^\infty_t L^2_x,  L^2_{x t})$. However, in general it is not known if such solutions stay bounded or regular for all $t>0$ for regular initial values. Neither the uniqueness of such solutions is known so far. A typical existence theorem for regular solutions always involves a small parameter in the initial condition, as a perturbation of a known regular solution.

Given a semi-open  time interval $[0, T)$, in this short note, we construct a regular solution of \eqref{nse} in a finite cylinder with a forcing term in the super critical space $L^{q}_t L^1_x$ for all $q>1$. The velocity is in the energy space at the final time $T$ when it blows up. Solutions beyond the energy space or super critical forcing terms have been investigated before. For a forcing term in the supercritical space $L^1_t L^2_x$, in the paper \cite{ABC22},   non-uniqueness of \eqref{nse} in the energy space is established. Since $2<\frac{3}{1}+\frac{2}{q}<\frac{3}{2} + \frac{2}{1}$ for large $q$, the forcing term here is actually milder than that in \cite{ABC22} with respect to the standard scaling.  Let us recall that a forcing term in local $L^q_tL^p_x$ space is super critical, critical and sub critical if $\frac{3}{p}+ \frac{2}{q} >, =, < 2$ respectively. Essentially speaking, under the standard scaling, a forcing term  is super critical, critical and sub critical if its scaling order is less than, equal to or greater to $-2$ respectively. In the homogeneous case $f=0$, non-uniqueness was proven for some solutions in the space $L^\infty_t L^2_x$ in \cite{BV}. Let us mention that an instantaneous  finite time blow up example was constructed in \cite{LPYZZZ} for a special cusp domain. The solution is in the energy space and the forcing term (Ampere force) is actually subcritical. These results indicate that if pertinent function spaces or the domain are sufficiently bad, then something bad can happen to the solutions.

Now let us describe the blow up solution in detail.
The solution is axially symmetric, namely  $v$ and $P$ are independent of the angle in a cylindrical coordinate system $(r,\,\th,\,x_3)$. That is, for $x=(x_1,\,x_2,\,x_3) \in \bR^3$,
$
r=\sqrt{x_1^2+x_2^2}, \quad
\th=\arctan (x_2/x_1),
$ and the basis vectors $e_r,e_\th,e_3$ are:
\[
e_r=(x_1/r, x_2/r, 0),\quad e_\th=(-x_2/r, x_1/r, 0),\quad e_3=(0,0,1).
\]
 In this case, solutions can be written in the form of
\be
\lab{vvrth3}
v=v_r(r,x_3, t)e_r+v_{\th}(r, x_3, t)e_{\th}+v_3(r, x_3, t)e_3
\ee and the forcing term can be written as
\be
\lab{ffrth3}
f=f_r(r,x_3, t)e_r+f_{\th}(r, x_3, t)e_{\th}+f_3(r, x_3, t)e_3
\ee

Therefore, $v_r$, $v_3$ and $v_\theta$ satisfy the axially symmetric Navier-Stokes equations
\be
\begin{aligned}
\lab{eqasns}
\begin{cases}
   \big (\Delta-\frac{1}{r^2} \big )
v_r-(v_r \p_r + v_3 \p_{x_3})v_r+\frac{(v_{\theta})^2}{r}-\partial_r
P-\p_t  v_r=f_r,\\
   \big   (\Delta-\frac{1}{r^2}  \big
)v_{\theta}-(v_r \p_r + v_3 \p_{x_3} )v_{\theta}-\frac{v_{\theta} v_r}{r}-
\partial_t v_{\theta}=f_{\theta},\\
 \Delta v_3-(v_r \p_r + v_3 \p_{x_3})v_3-\p_{x_3} P-\p_t v_3=f_3,\\
 \frac{1}{r} \p_r (rv_r) +\p_{x_3}
v_3=0,
\end{cases}
\end{aligned}
\ee which will be abbreviated as ASNS.

If the swirl $v_\theta=0$ and $f$ is sufficiently regular, then it is well known (
\cite{L},  \cite{UY}), that finite energy solutions
to the Cauchy problem of (\ref{eqasns}) in $\mathbb{R}^3$ are smooth for all time $t>0$.
See also \cite{LMNP}.
 Some existence results with a small parameter can be found in \cite{HL}, \cite{HLL} and \cite{ZZ}.  Critical or sub-critical regularity conditions can be found in \cite{NP}, \cite{CL} and \cite{JX}.
Despite many difficulty, there has been no lack of research efforts on ASNS. Here is a partial list of recent papers on conditional regularity of ASNS: \cite{CSTY1},
\cite{CSTY2},   \cite{KNSS}, \cite{LZ11},
\cite{SeZh},
 \cite{Panx}, \cite{CFZ}, \cite{LZ17}, \cite{Weid}.
Comparing with the full 3 dimensional NS, it was realized in \cite{LZ17} that
the vortex stretching term of the ASNS is critical after a  suitable change of dependent variables.

Due to the presence of a boundary, we must specify suitable boundary values.
Let $D$ be a piecewise smooth domain in $\mathbb{R}^3$ and $n$ be the unit outward normal of the smooth part of $\p D$. Recall that the Navier \cite{Nav} (total) slip boundary condition reads
\be
\lab{nvslbc}
(\mathbb{S}(v) n)_{tan} = 0, \qquad v \cdot n = 0, \quad \text{on} \quad \p D.
\ee  Here $\mathbb{S}(v)=\frac{1}{2} (\na v + (\na v)^T)$ is the strain tensor and $(\mathbb{S}(v) n)_{tan}$ stands for the tangential component of the vector $\mathbb{S}(v) n$. In the axially symmetric setting, the matrix
\be
\lab{Sv}
 \mathbb{S}(v) = \begin{bmatrix} \p_r v_r & \frac{1}{2} (\p_r v_\th - \frac{1}{r} v_\th) & \frac{1}{2}(\p_{x_3} v_r + \p_r v_3) \\
                          \frac{1}{2} (\p_r v_\th - \frac{1}{r} v_\th) &  \frac{1}{r} v_r & \frac{1}{2} \p_{x_3} v_\th \\
                           \frac{1}{2}(\p_{x_3} v_r + \p_r v_3) & \frac{1}{2} \p_{x_3} v_\th &  \p_{x_3} v_3
\end{bmatrix}
\ee is the strain tensor. See \cite{Zq22} Section 1, for a quick derivation e.g.

If $D$ is a finite cylinder with $x_3$ axis as its axis of rotation, then the boundary
$\p D$ can be written as the union of horizontal and vertical parts, which are denoted by $\p^H D$ and $\p^V D$ respectively.  From  \eqref{Sv} and \eqref{nvslbc}, one sees that the Navier slip boundary condition can be expressed explicitly as
\be
\lab{nvbc2}
\bali
\p_{x_3} v_r=\p_{x_3} v_\th=0, \quad v_3=0, \quad \text{on} \quad \p^H D;\\
\p_r v_\th=\frac{v_\th}{r}, \, \p_r v_{x_3}=0, v_r=0, \quad \text{on} \quad \p^V D.
\eali
\ee

The main result of the note is part (2) of the following

\begin{theorem}
\lab{prblows}
Let $D$ be the cylinder $B_2(0, 1) \times [0, 1]$ in $\mathbb{R}^3$, where $B_2(0, 1)$ is the unit ball in $\mathbb{R}^2$ centered at the origin.

(1). Let $k=k(r)$ be a smooth nontrivial function supported in the unit interval $[0, 1]$ and $\alpha=\int^1_0 s e^{\frac{s^2}{2}} \int^\infty_s e^{- \frac{1}{2} l^2} k(l) dl ds$. Given $T \in (0, 1/2]$, the vector field
 \[
 v=v(x, t)= -\frac{1}{r} \int^{r/\sqrt{2(T-t)}}_0 s e^{\frac{s^2}{2}} \int^\infty_s e^{- \frac{1}{2} l^2} k(l) dl ds \, e_\theta + \alpha r \, e_\theta
 \]is an unbounded classical solution of the forced axially symmetric Navier Stokes equation:
\be
\lab{nsF}
\Delta v -v \nabla v -\nabla P -\p_t v =\frac{1}{(2(T-t))^{3/2}} k\left(\frac{r}{\sqrt{2(T-t)}} \right)  e_\th \quad \text{in} \quad D \times [0, T),
\ee which satisfies the no slip boundary condition $v=0$ on the vertical boundary of $D$ and the Navier total slip boundary condition on the horizontal boundary of $D$. In addition $v$ satisfies
 \[
\int_{D} |v|^2(x, t) dx +
\int^t_0 \int_{D} |\na v|^2(x, s) dx ds \le C |\ln (T-t)|,
\]and the forcing term is in the space $L^{2^-}_t L^1_x$.

(2). In addition, suppose $k=k(r) \le 0$ and $k \neq 0$. Then
\be
\overline{v} =\left[\ln \left(-\frac{1}{r} \int^{r/\sqrt{2(T-t)}}_0 s e^{\frac{s^2}{2}} \int^\infty_s e^{- \frac{1}{2} l^2} k(l) dl ds+1 \right) - ( \ln (1-\alpha)) \, r \right] e_\theta
\ee is a classical unbounded solution to
\[
\Delta v -v \nabla v -\nabla P -\p_t v =Y  e_\th \quad \text{in} \quad D \times [0, T)
\]with some $Y \in L^{q}_{t} L^1_x, \quad \forall q>1$.
Moreover $\overline{v}$ is in the energy space up to $T$ and satisfies the Navier slip boundary condition on the horizontal boundary of $D$ and the no slip condition on the vertical boundary of $D$.
\proof
\end{theorem}

(1).
This part, as a preparation and motivation for part (2),  can be checked directly but we prefer to show how the solution is found.
After setting $v_r=0, v_3=0, f_r=0, f_3=0, f_\theta=f_\theta(r, t)$ in the ASNS \eqref{eqasns} with $v$ and $f$ given \eqref{vvrth3} and \eqref{ffrth3},  it is reduced to:
\be
\begin{aligned}
\lab{pswasns}
\begin{cases}
   \frac{(v_{\theta})^2}{r}-\partial_r
P=0,\\
   \big   (\Delta-\frac{1}{r^2}  \big
)v_{\theta}-
\partial_t v_{\theta}=f_\theta,\\
\p_{x_3} P=0.
\end{cases}
\end{aligned}
\ee Let $\phi$ be a solution to the following equation such that $\phi(0, t)=0$.
\be
\lab{eqf}
\left(\Delta - \frac{1}{r^2}\right) \phi(r, t) -\p_t \phi(r, t)=h
\ee in $D \times [0, T]$, $T>0$.

Set $v_\th=\phi=\phi(r, t)$, $h=f_\theta(r, t)$ and
\[
P=\int^r_0 \frac{v^2_\th(l, t)}{l} dl,
\]which is well defined  since $v_\th(0, t)=0$ when $v$ is smooth. One checks easily that the pair $v, P$ satisfy \eqref{pswasns}. Therefore $v=v_\th \, e_\th$ satisfies ASNS.

We further look for self similar solutions of \eqref{eqasns} in the form of
\be
\lab{sss}
\al
v&=v(x, t) = \frac{1}{\sqrt{2(T-t)}} V\left(\frac{r}{\sqrt{2(T-t)}} \right),  \qquad f=f(x, t) = \frac{1}{[2(T-t)]^{3/2}} F\left(\frac{r}{\sqrt{2(T-t)}} \right),\\
P&=P(x, t) = \frac{1}{2(T-t)} P_0\left(\frac{r}{\sqrt{2(T-t)}} \right).
\eal
\ee

Then solutions $\phi$ to \eqref{eqf} is reduced to $\phi=\frac{1}{\sqrt{2(T-t)}} \phi_0\left(\frac{r}{\sqrt{2(T-t)}} \right)$ where $\phi_0$ solves the ode:
\be
\lab{fode}
\phi''_0+ \frac{1}{r} \phi'_0 - \frac{1}{r^2} \phi_0 -  \phi_0 -  r \phi'_0=k(r), \qquad \phi_0(0)=0, \quad k(r)=F \cdot e_\theta, \qquad r \ge 0.
\ee Write $g= r \phi_0$. By direct calculation, $g$ satisfies the equation
\[
g'' -  \left( \frac{1}{r} + r \right) g' = r k(r), \qquad g'(0)=0.
\]Therefore
\[
\left( e^{- \ln r - \frac{1}{2} r^2} g'(r) \right)' = e^{- \ln r - \frac{1}{2} r^2} r k(r),
\]which shows
\[
g(r)= - \int^r_0 s e^{\frac{s^2}{2}} \int^\infty_s e^{- \frac{1}{2} l^2} k(l) dl ds.
\]i.e.
\be
\lab{formf}
\phi_0=\phi_0(r) = -\frac{1}{r} \int^r_0 s e^{\frac{s^2}{2}} \int^\infty_s e^{- \frac{1}{2} l^2} k(l) dl ds.
\ee So a solution $u$ to \eqref{nsF} is given explicitly by
\be
\lab{u=phietc}
u=\phi e_\theta=\frac{1}{\sqrt{2(T-t)}} \phi_0\left(\frac{x}{\sqrt{2(T-t)}} \right)  e_\theta
=-\frac{1}{r} \int^{r/\sqrt{2(T-t)}}_0 s e^{\frac{s^2}{2}} \int^\infty_s e^{- \frac{1}{2} l^2} k(l) dl ds \, e_\theta.
\ee

Since $k=k(r)$ is smooth and compactly supported, it is clear that there is a positive constant $C$ depending only on $k$ such that
\[
|\phi_0(r)| \le \frac{C r}{r^2 + 1}, \qquad |\phi'_0(r)| \le \frac{C }{r^2 + 1}.
\]Therefore
\be
\lab{sssu}
|u(r, t)| \le \frac{C r}{r^2 + T-t}, \qquad |\nabla u(r, t)| \le \frac{C }{r^2 + T-t}.
\ee Here and later $C$ may change by a numerical factor from line to line. From these,
 one can easily deduce
\be
\lab{eneru}
\int_{D} |u|^2(x, t) dx +
\int^t_0 \int_{D} |\na u|^2(x, s) dx ds \le C |\ln (T-t)|.
\ee Thus $u=u(x, t)$ is in the energy space in $D \times [0, T^{-}]$ and logarithmically beyond the energy space in  $D \times [0, T^{-}]$.

Consider the vector field
\be
\lab{v=u+a}
v \equiv u +\alpha r e_\theta.
\ee Observe that $r$ is a stationary solution to $(\Delta-\frac{1}{r^2}-\partial_t) r =0$.
 Hence $\phi+ \alpha r$ is also a solution of \eqref{pswasns}, which shows $v= (\phi+ \alpha r) e_\theta$ a solution of the Navier Stokes equations \eqref{nsF} with suitably adjusted pressure $P$. Due to \eqref{eneru}, the following holds
 \be
\lab{enerv}
\int_{D} |v|^2(x, t) dx +
\int^t_0 \int_{D} |\na v|^2(x, s) dx ds \le C |\ln (T-t)|.
\ee

 Since $v$ is independent of the vertical variable $x_3$ and $v_r=v_3=0$, we see that $v$ satisfies the Navier total slip boundary condition on  the horizontal boundary of $D$.  On the vertical boundary $r=1$, the velocity $v$ is $0$ since $u=-\alpha e_\theta$. The reason is that $0\le T-t \le 1/2$ by the choice of $T$ and,  if $r=1$, then $r/\sqrt{2(T-t)} \ge 1$. Since the function $k$ is supported in $[0, 1]$, we have
 \be
 \lab{u1t}
 \al
 u(1, t) &= \phi(1, t) e_\theta = -\int^{r/\sqrt{2(T-t)}}_0 s e^{\frac{s^2}{2}} \int^\infty_s e^{- \frac{1}{2} l^2} k(l) dl ds \, e_\theta\\
 &= -\int^1_0 s e^{\frac{s^2}{2}} \int^\infty_s e^{- \frac{1}{2} l^2} k(l) dl ds \, e_\theta = -\alpha e_\theta.
 \eal
\ee Therefore the no slip boundary condition is satisfied on the vertical boundary.

Finally, since $\int_D \frac{1}{(T-t)^{3/2}} |k|\left(\frac{r}{\sqrt{2(T-t)}}\right)  dx \le \frac{C}{\sqrt{T-t}}$, we see that the forcing term is in $L^{2^-}_t L^1_x$. This proves (1).

(2). We write $\eta= \ln (\phi+1)$ where $\phi$ is given in \eqref{u=phietc}. Then the following identity holds
\be
\lab{eqeta}
\al
\left(\Delta - \frac{1}{r^2}\right) \eta(r, t) -\p_t \eta(r, t)
&=-\frac{1}{r^2} \eta + \frac{\phi}{r^2 (\phi+1)} + \frac{h}{\phi+1} - \left( \frac{ \partial_r \phi}{\phi+1} \right)^2 \\
&\equiv Y \equiv Y_1+Y_2+Y_3+Y_4,
\eal
\ee
where $h=\frac{1}{[2(T-t)]^{3/2}} k\left(\frac{r}{\sqrt{2(T-t)}} \right)$.

Since $\phi(1, t) = -\alpha>0$ from \eqref{u1t}, we have $\eta(1, t) = \ln (1-\alpha)$. Therefore, as in part (1), the vector field
\be
\overline{v} = \left(\eta(r, t) - [\ln (1-\alpha)] \, r \right) e_\theta=\left[\ln (\phi(r, t)+1) - (\ln (1-\alpha)) \, r \right] e_\theta
\ee is a classical  solution to
\[
\Delta v -v \nabla v -\nabla P -\p_t v =Y  e_\th \quad \text{in} \quad D \times [0, T),
\] which blows up at time $T$. It also satisfies the Navier slip boundary condition on the horizontal boundary of $D$ and the no slip condition on the vertical boundary of $D$.

Now let us check that $\overline{v}$ is in the energy space:
\be
\lab{engvbar}
\int_{D} |\overline{v}|^2(x, t) dx +
\int^t_0 \int_{D} |\na \overline{v}|^2(x, s) dx ds \le C, \quad  t \in [0, T].
\ee Using \eqref{sssu}, we have, by $\ln^2(1+a) \le C a$ for $a>0$, that
\[
\int_D |\overline{v}|^2(x, t) dx \le C \int^1_0 \ln^2\left(1+ \frac{r}{r^2+T-t} \right) r dr +C \le  C \int^1_0 \frac{r^2}{r^2+T-t}  dr +C \le 2C;
\]
\[
\al
\int^T_0 \int_{D} &|\na \overline{v}|^2(x, t) dx dt \le \int^T_0 \int_{D} \left( \frac{ \partial_r \phi}{\phi+1} \right)^2 dx dt + C \int^T_0 \int_{D} \frac{1}{r^2} \ln^2\left(1+ \frac{r}{r^2+T-t} \right) dxdt+C  \\
&\le C \int^T_0 \int^1_0 (r^2+T-t)^{-2}\left(\frac{r}{r^2+T-t} + 1 \right)^{-2} r drdt+C
\\
&=C \int^T_0 \int^1_0 \frac{r}{(r+r^2+T-t)^2}   dr dt+C \le C \int^1_0 \frac{r}{r+r^2} dr +C \le 2 C.
\eal
\]Here we have used the fact that
\be
\lab{phi>}
\phi(r, t)  \ge \frac{ r}{C(r^2 + T-t)}
\ee which is due to \eqref{formf} and \eqref{u=phietc} and the assumption $- k \ge 0$, $k \neq 0$ so that the integral $-\int^\infty_s e^{-l^2/2} k(l) dl$ converges to a positive constant as $ s \to 0$. This also shows that $\overline{v}$ blows up as $ t \to T$.

Therefore the energy inequality \eqref{engvbar} is valid.

Next, by \eqref{sssu} again
\[
\al
&\int_D |Y_1| dx  \le C \int^1_0 \frac{1}{r^2} \ln\left(1+ \frac{C r}{r^2+T-t} \right) r dr \qquad \qquad (l \equiv r/\sqrt{T-t})\\
&=C \int^{1/\sqrt{s}}_0 \frac{1}{l} \ln \left(1+\frac{1}{\sqrt{s}} \frac{l}{l^2+1} \right)dl=C \int^{\sqrt{s}}_0 ... dl  + C \int^{1/\sqrt{s}}_{\sqrt{s}} ... dl
\qquad \qquad (s=T-t)\\
& \le C \int^{\sqrt{s}}_0 \frac{1}{\sqrt{s}(l^2+1)} dl + C \int^{1/\sqrt{s}}_{\sqrt{s}}
\frac{1}{l} \ln \left(1+\frac{1}{\sqrt{s}}\right) dl
\le C \ln^2(T-t).
\eal
\]This shows $Y_1 \in L^{q}_t L^1_x$ for all $q>1$. Similarly
\[
\int_D |Y_2| dx \le  C \int^1_0 \frac{1}{r+r^2+T-t} dr \le C \ln\frac{1}{T-t};
\]
\[
\al
\int_D |Y_3| dx &\le C \int^1_0 \frac{r^2+T-t}{r+r^2+T-t} \frac{1}{[2(T-t)]^{3/2}} k\left(\frac{r}{\sqrt{2(T-t)}} \right) r dr, \qquad (l \equiv r/\sqrt{2(T-t)}) \\
&\le \frac{C }{\sqrt{T-t}} \int^1_0 \frac{l^2+1}{(l/\sqrt{T-t}) + l^2} k(l) l dl
 \le C, \qquad ( supp \,  k \subset [0,  1]);
\eal
\]
\[
\int_D |Y_4| dx \le C \int^1_0 \frac{r}{(r+r^2+T-t)^2}   dr \le C\int^1_0 \frac{r}{r^2+(T-t)^2}   dr \le C \ln\frac{1}{T-t}.
\]Therefore $Y=Y_1+Y_2+Y_3+Y_4 \in L^q_t L^1_x$ for all $q>1$.  By direct observation, the scaling order of the forcing term $Y$  is slightly worse than $-2$, the critical order.

The proof of the theorem is complete.

\qed

We remark that if one allows arbitrary forcing terms, then one can produce blow up solutions with slightly super critical forcing term in the localized space $L^\infty_t L^{3^-/2}_x$. Here $3^-$ is any positive number strictly less than $3$. Recall that for the forcing term, the space $L^\infty_t L^{3/2}_x$ is critical. The following is an example. Let $\phi=\phi(y)$ be a smooth, divergence free vector field, which is compactly supported in $D$. For example $\phi = \text{curl} \, X$ where $X$ is any smooth vector field compactly supported in $D$.
Define
\[
v=v(x, t) = \ln (1/(1-t)) \phi (x/\sqrt{1-t}).
\]Then $v$ is trivially a blow up solution of the Navier-Stokes equation with $0$ pressure and $0$ boundary value in $D \times [0, 1)$.
\[
\Delta v - v \nabla v - \partial_t v =F
\]where $F=\Delta v - v \nabla v - \partial_t v$. It is easy to check that $v$ is in the energy space and $F \in L^\infty_t L^{3^-/2}_x([0, 1] \times D)$. This example is inspired by an example by an anonymous reviewer of a journal who suggested the self similar solution $u=(1/\sqrt{1-t}) \phi (x/\sqrt{1-t})$, then the corresponding forcing term will be in the super critical space $L^\infty_t L^1_x$. However unlike the rotational forcing terms in this note, these forces are unlikely to be physically achievable.

\section*{Acknowledgments} We wish to thank Prof. Hongjie Dong,  Zhen Lei, Zijin Li,
  Xinghong Pan, Xin Yang and Na Zhao and Dr. Chulan Zeng for helpful discussions. The support of Simons Foundation grant 710364 is gratefully acknowledged.

\bibliographystyle{plain}



\def\cprime{$'$}

\end{document}